\documentclass[11pt]{amsart}

\usepackage{pb-diagram}

\newtheorem{thm}{Theorem}[section]    

\newtheorem{lem}[thm]{Lemma}
\newtheorem{prop}[thm]{Proposition}

\theoremstyle{remark}                
\newtheorem{rmk}[thm]{Remark}

\theoremstyle{definition}

\numberwithin{equation}{section}

\newcommand{\secref}[1]{Section~\ref{#1}}
\newcommand{\thmref}[1]{Theorem~\ref{#1}}

\newcommand{\propref}[1]{Proposition~\ref{#1}}
\newcommand{\lemref}[1]{Lemma~\ref{#1}}

\newcommand{\clsp}{\overline{\operatorname{span}}}

\newcommand{\id}{\operatorname{id}}
\newcommand{\spn}{\operatorname{span}}
\newcommand{\Ad}{\operatorname{Ad}}

\newcommand{\slashstyle}{\mathcal}

\renewcommand{\L}{{\slashstyle L}}

\newcommand{\K}{{\slashstyle K}}

\newcommand{\FF}{\mathcal F}

\newcommand{\bbC}{\mathbb C}

\newcommand{\midtext}[1]{\quad\text{#1}\quad}

\newcommand{\su}{S_U}
\newcommand{\sv}{S_V}
\newcommand{\sw}{S_W}

\newcommand{\intextmatrix}[1]
 {\left(\begin{smallmatrix}#1\end{smallmatrix}\right)}
\newcommand{\imp}[2]{$#1$ -- $#2$}

\newcommand{\inner}[1]{\left\langle #1 \right\rangle}
\newcommand{\rinner}[2]{\left\langle #1 \right\rangle_{#2}}
\newcommand{\linner}[2]{\,{}_{#2}\!\left\langle #1 \right\rangle}

\newcommand{\twoby}[4]{\begin{pmatrix}#1&#2\\#3&#4\end{pmatrix}}
\newcommand{\texttwoby}[4]{\left(\begin{smallmatrix}#1&#2\\
#3&#4\end{smallmatrix}\right)}
\newcommand{\ttwoby}{\texttwoby}

\begin{document}

\title[Equivariance and Imprimitivity]{Equivariance and Imprimitivity 
for Discrete Hopf $C^*$-Coactions}

\date{August 10, 1998.}

\author[Kaliszewski]{S.~Kaliszewski}
\address{Department of Mathematics \\ Arizona State University \\
Tempe, AZ 85287}
\email{kaz@math.la.asu.edu}

\author[Quigg]{John Quigg}
\address{Department of Mathematics \\ Arizona State University \\
Tempe, AZ 85287}
\email{quigg@math.la.asu.edu}

\subjclass{46L55, secondary 22D25}

\thanks{This research was partially supported by the Australian
Research Council and by the National Science Foundation (under Grant
No. DMS9401253).}

\begin{abstract}
Let $U$, $V$, and $W$ be multiplicative unitaries coming from discrete
Kac systems such that $W$ is an amenable normal submultiplicative
unitary of $V$ with quotient $U$.  We define notions for right-Hilbert
bimodules of coactions of $\sv$ and $\hat\sv$, their restrictions to $\sw$
and $\hat\su$, their dual coactions, and their full and reduced crossed
products.  If $N(A)$ denotes the imprimitivity bimodule associated to
a coaction $\delta$ of $\sv$ on a $C^*$-algebra $A$ by Ng's
imprimitivity theorem, we prove that 
for a suitably nondegenerate injective
right-Hilbert bimodule coaction of $\sv$ on $_AX_B$, 
the balanced tensor products 
$N(A)\otimes_{A\times\hat\sw}(_AX_B\times\hat\sw)$ and 
$(_AX_B\times\hat\sv\times_r\su)\otimes_{B\times\hat\sv\times_r
\su}N(B)$
are isomorphic 
right-Hilbert $A\times\hat\sv\times_r\su$ -- $B\times\hat\sw$
bimodules.  This can be 
interpreted as a natural equivalence between certain
crossed-product functors.
\end{abstract}

\maketitle


\section{Introduction}
\label{intro-sec}

Since Baaj and Skandalis introduced multiplicative unitaries in
\cite{BS-UM} as a generalization of locally compact groups, and proved a
duality theorem (\cite[Th\'eor\`eme~7.5]{BS-UM})
for crossed products by
coactions of the associated Hopf $C^*$-algebras, there has been much
interest in extending other results for group actions and coactions to
this context.  Recently Ng (\cite{NgME}) has defined notions of sub- and
quotient multiplicative unitaries, and has proved that for
multiplicative unitaries $U$, $V$, and $W$ coming from discrete
Kac systems such that
$W$ is an amenable normal submultiplicative 
unitary of $V$ with
quotient $U$, and for any injective nondegenerate coaction $\delta$ of
$\sv$ on a $C^*$-algebra $A$, the iterated crossed product
$A\times_\delta\hat{S}_V\times_{\hat\delta|,r} S_U$ is Morita
equivalent to $A\times_{\delta|}\hat{S}_W$ (\cite[Theorem~3.4]{NgME}).  
This is an analog both of Green's celebrated imprimitivity theorem
(\cite{GreLS}), which implies that for an action $\alpha$ of a group $G$
on $A$ and a closed normal subgroup $N$ of $G$, $A\times_\alpha
G\times_{\hat\alpha|}G/N$ is Morita equivalent to $A\times_{\alpha|}N$,
and of Mansfield's imprimitivity theorem (\cite{ManIR2}) for coactions (as
generalized to non-amenable groups in \cite{KQ-IC}), which provides a
Morita equivalence between $A\times_\delta G\times_{\hat\delta|,r}N$ and
$A\times_{\delta|}G/N$ for a coaction $\delta$ (satisfying a mild
condition)
of $G$ on $A$ and any
closed normal subgroup $N$ of $G$.  
For discrete multiplicative unitaries, Ng's theorem generalizes
Baaj-Skandalis duality (ignoring differences between full and reduced
crossed products) in the same way that Green's theorem generalizes
Imai-Takai-Takesaki duality (\cite{IT-DC}), and Mansfield's theorem
generalizes the duality of Katayama (\cite{KatTD}).

Now the significance of Green's theorem is that his imprimitivity
bimodule may be viewed as a Hilbert $A\times_{\alpha|} N$-module with a
nondegenerate left action of $A\times_\alpha G$
by adjointable operators, and thus allows induction of
representations from $A\times_{\alpha|} N$ to $A\times_\alpha G$ via 
Rieffel's framework
(\cite{RieIR}).  Similarly, Mansfield's bimodule allows induction of
representations from $A\times_{\delta|}G/N$ to $A\times_\delta G$.  
The representation-inducing processes arising from these bimodules, and
their interactions with one another, have received much
attention lately (see
\cite{EchDI, GL-AN, EKR-CP, KQR-DR, NilDC}),
and the method that has evolved is to work with the bimodules that
implement
the inducing maps on representations,
rather than with those inducing
maps themselves.  We call the bimodules involved {\em right-Hilbert
bimodules\/}; they are essentially imprimitivity bimodules ${}_KX_B$
together with
nondegenerate homomorphisms of $A$ into $M(K)$.

The equivariant right-Hilbert bimodules --- that is, those
right-Hilbert bimodules
$_AX_B$ which carry compatible actions or coactions of a group $G$ ---
turn out to be closely related to imprimitivity theorems.  In work
with S.~Echterhoff and I.~Raeburn which is currently in preparation
we have shown, for example, that Green's imprimitivity theorem can be
viewed as a natural equivalence
between the crossed product functors
$(A,G,\alpha)\mapsto A\times_\alpha G\times_{\hat\alpha|}G/N$ and
$(A,G,\alpha)\mapsto A\times_{\alpha|}N$ defined on a category whose
objects are $C^*$-algebras with actions of $G$ and whose morphisms
$(A,\alpha)\to (B,\beta)$ are (isomorphism classes) of equivariant
right-Hilbert $A$ -- $B$ bimodules (\cite{EKQR}).  

In this paper, we show that Ng's imprimitivity theorem is similarly
compatible with equivariant right-Hilbert bimodules.  To do so, we must
first develop a theory of coactions of Hopf $C^*$-algebras $\sv$ and
$\hat\sv$ on right-Hilbert bimodules, and their crossed products; this
is done as efficiently as possible in \secref{rh-bimod-sec} by building for
the most part on Ng's imprimitivity bimodule apparatus (\cite{NgCC2}).
In \secref{fixed-point}, we review Ng's fixed-point theorem
(\cite[Proposition~2.11]{NgME}), since it provides the construction of the
bimodule which appears in his
imprimitivity theorem.  Here we prove two
lemmas relating Ng's bimodule to the linking algebra and standard
right-Hilbert bimodule (see below) constructions we use in proving our
main theorem.  

In the final section, we prove our main result: for $U$,
$V$, and $W$ as in Ng's theorem, and for a suitably nondegenerate
injective right-Hilbert bimodule coaction of $\sv$ on $_AX_B$,
$$N(A)\otimes_{A\times\hat{S}_W}X\times\hat{S}_W \cong
X\times\hat{S}_V\times_r S_U\otimes_{B\times\hat{S}_V\times_r S_U} N(B)$$
as right-Hilbert
$A\times\hat{S}_V\times_r S_U$ -- $B\times\hat{S}_W$ bimodules, where
$N(A)$ denotes Ng's $A\times\hat S_V\times_r S_U$ -- $A\times\hat
S_W$ imprimitivity bimodule (and similarly for $N(B)$). 
As discussed above for group actions, this should give
a natural equivalence between certain crossed-product functors, although
we don't formalize this in the present paper. 
(Part of our point
here is that any reasonable imprimitivity theorem should be compatible
with equivariant right-Hilbert bimodules, and that the proof of this,
following the same strategy we use in the proof of \thmref{main}, should
be relatively straightforward.)  Our theorem should 
have implications for induced and restricted representations of
crossed products by Hopf $C^*$-algebras, and for equivariant $KK$-theory
as in \cite{BS-CH}.

A substantial amount of this research was conducted while the first
author visited the second at Arizona State University, and on another
occasion while both authors visited the University of Newcastle. We
particularly thank Iain Raeburn in Newcastle for his hospitality.

\section*{Preliminaries}

For compatibility with Ng's work on imprimitivity bimodules, we define
a {\em right-Hilbert $A$ -- $B$ bimodule\/} over $C^*$-algebras $A$
and $B$ to be an imprimitivity bimodule ${}_KX_B$ together with a
nondegenerate homomorphism of $A$ into $M(K)$.  If $X$ is a full
Hilbert $B$-module (cf. \cite{LanHC}), then $X$ is a $\K_B(X)$ -- $B$
imprimitivity bimodule, and $M(\K_B(X))=\L_B(X)$, so this is the same as
having a nondegenerate action of $A$ by adjointable operators on $X$.
Note that $K$ 
itself becomes a right-Hilbert $A$ -- $K$ bimodule by using the
natural $K$ -- $K$ imprimitivity bimodule structure on $K$; we call this
a {\em standard\/} right-Hilbert bimodule.  We have the decomposition
${}_AX_B\cong {}_AK\otimes_KX_B$ of any right-Hilbert bimodule as a
balanced tensor product of a standard bimodule and an imprimitivity
bimodule. We use the conventions
of \cite{ER-MI} 
regarding multiplier bimodules, linking algebras,
and homomorphisms of imprimitivity bimodules.

Let $(S,\delta_S)$ be a Hopf $C^*$-algebra, and let $\delta\colon
A\to M(A\otimes S)$ be a coaction of $S$ on a $C^*$-algebra $A$, as in
\cite[Definition 0.2]{BS-UM}. The coaction
$\delta$ is called {\em
nondegenerate\/} if $\clsp\{\delta_A(A)(1\otimes S)\} = A\otimes S$.
A {\em covariant pair\/} for $(A,S,\delta)$ on a $C^*$-algebra $B$
consists of a nondegenerate homomorphism $\theta\colon A\to M(B)$ and
a unitary corepresentation $u\in M(B\otimes S)$ of $S$ 
such that
$$(\theta\otimes\id)\circ\delta(a) = \Ad(u)(\theta(a)\otimes 1)$$
for each $a\in A$ 
(\cite[Definition~2.8]{NgCC1}).  The {\em full crossed product\/} for
$(A,S,\delta)$ is a $C^*$-algebra $A\times_\delta\hat S$ together with
a universal covariant pair $(j,v)$ for $(A,S,\delta)$ on
$A\times_\delta\hat S$ (\cite[Definition~2.11(b)]{NgCC1}).
If
$S=\hat\sv$ for a multiplicative unitary $V$ coming from a Kac system
(see below), we write $A\times_\delta\sv$ 
(with no hat) for $A\times_\delta\hat S$.  

Let $V\in \L(H\otimes H)$ be a regular multiplicative unitary 
as in \cite{BS-UM}.  
We let $L$ and $\rho$ denote the maps of $\L(H)_*$ into $\L(H)$ defined
by 
$$L(\omega) = (\omega\otimes\id)(V)\midtext{and} \rho(\omega) =
(\id\otimes\omega)(V);$$
then we have the associated reduced Hopf $C^*$-algebras 
$$\sv = \clsp\{ L(\omega)\mid \omega\in\L(H)_*\} \midtext{and}
\hat\sv = \clsp\{ \rho(\omega)\mid \omega\in\L(H)_*\}$$
with comultiplications $\delta_V$ and $\hat\delta_V$, respectively,
given by 
$$\delta_V(x) = V(x\otimes 1)V^*\midtext{and} 
\hat\delta_V(y) = V^*(1\otimes y)V$$
(\cite[Th\'eor\`eme~3.8]{BS-UM}).  
The corresponding
full Hopf $C^*$-algebras are denoted $(S_V)_p$ and $(\hat\sv)_p$, but
their comultiplications are still denoted
$\delta_{V}$ and $\hat\delta_{V}$ (\cite[Corollaire~A.6]{BS-UM}). 
We view $L$ both as a faithful representation of $\sv$ and as a
nondegenerate representation of $(\sv)_p$ on $\L(H)$; similarly, we
view $\rho$ as a map on both $\hat\sv$ and $(\hat\sv)_p$ (cf.
\cite[Proposition~1.16(i)]{NgCC1}).  

The unitary corepresentations $u\in M(B\otimes\sv)$ of $\sv$ are in
bijective correspondence with the nondegenerate homomorphisms
$\nu\colon (\hat\sv)_p\to M(B)$ (\cite[Lemma2.6]{NgCC1}).  
If $(A,\sv,\delta)$ is a coaction, by \cite[Remark~2.12(b)]{NgCC1} we have 
$$A\times_\delta\hat\sv = \clsp\{ j(a)\mu(y) \mid a\in A, y\in
(\hat\sv)_p \},$$
where $\mu\colon (\hat\sv)_p\to M(A\times_\delta\hat\sv)$ is the
nondegenerate homomorphism corresponding to~$v$. 
For every pair of homomorphisms $\theta\colon A\to M(B)$ and
$\nu\colon(\hat\sv)_p\to M(B)$ coming from a covariant pair $(\theta,u)$,
there is (by definition of the crossed product)
a unique nondegenerate homomorphism $\theta\times\nu\colon
A\times_\delta\hat\sv\to M(B)$ such that
$$(\theta\times\nu)\circ j = \theta \midtext{and}
(\theta\times\nu)\circ \mu = \nu,$$
and the latter condition is equivalent to
$$((\theta\times\nu)\otimes\id)(v) = u.$$

Let
$\pi_L = (\id\otimes L)\circ\delta\colon A\to
\L_A(A\otimes H)$. Then the {\em reduced crossed
product\/} (\cite[Definition~2.11(a)]{NgCC1}) is
$$A\times_{\delta,r}\hat\sv = C^*\bigl(\{ \pi_L(a)
(1\otimes\rho(\omega)) \mid a\in A, \omega\in \L(H)_*
\}\bigr) \subseteq \L_A(A\otimes H).$$
The reduced crossed product by a coaction
$\delta_p$ of $(\sv)_p$ is defined similarly, and we have
$$A\times_{\delta_p,r}(\hat\sv)_p\cong A\times_{\delta,r}\hat\sv,$$
where
the coaction $\delta=(\id\otimes L)\circ\delta_p$ 
of $\sv$ is the \emph{reduction} of $\delta_p$
(\cite[Proposition~2.14]{NgCC1}).
There is
a {\em dual coaction\/} $\hat\delta$ of
$(\hat\sv)_p$ on the full crossed
product $A\times_\delta\hat\sv$ which satisfies
$$\hat\delta(j(a)\mu(y)) = (j(a)\otimes
1)(\mu\otimes\id)(\hat\delta_{V}(y))$$
for all $a\in A$, $y\in (\hat\sv)_p$; we also denote its reduction by
$\hat\delta$ (\cite[Proposition~2.13]{NgMM}).  

Now suppose $V$ comes from a Kac system $(H,V,U)$
(\cite[D\'efinition~6.4]{BS-UM}). 
For any coaction $\delta$ of $\sv$ on $A$, 
we have
\begin{equation}\label{xr-eq1}
A\times_{\delta,r}\hat\sv = \clsp\{ \pi_L(a)
(1\otimes\rho(\omega)) \mid a\in A, \omega\in \L(H)_*
\} \subseteq \L_A(A\otimes H).
\end{equation}
Similarly, for any coaction $\delta$ of
$\hat\sv$ on $A$, we denote the reduced crossed product by
$A\times_{\delta,r}\sv$ (\cite[D\'efinition~7.1]{BS-UM}), and we have
\begin{equation}\label{xr-eq2}
A\times_{\delta,r}\sv = \clsp\{ \hat\pi_\lambda(a) 
(1\otimes L(\omega)) \mid a\in A, \omega\in \L(H)_*
\} \subseteq \L_A(A\otimes H),
\end{equation}
where $\lambda = \Ad(U)\circ\rho$ and
$\hat\pi_\lambda = (\id\otimes\lambda)\circ\delta$
(\cite[Lemme~7.2]{BS-UM}).  

It is important to note that for any Kac system $(H,V,U)$,
$(H,\hat V,U)$ is also a Kac system (\cite[Proposition~6.5]{BS-UM}),
where $\hat V = \Sigma (U\otimes 1)V(U\otimes 1)\Sigma$ and
$\Sigma$ denotes the flip operator on $H\otimes H$, and that then
$(\hat\sv,\hat\delta_V)\cong (S_{\hat V}, \delta_{\hat V})$ as a
consequence of \cite[Proposition~6.7]{BS-UM}.
Hence, for any coaction $(A,\hat\sv,\delta)$ there is
a coaction $(A,S_{\hat V},\delta')$ such that $A\times_\delta\sv\cong
A\times_{\delta'}\hat S_{\hat V}$ and $A\times_{\delta,r}\sv\cong
A\times_{\delta',r}\hat S_{\hat V}$.  Thus any results about crossed
products by coactions of $\sv$ always yield analogous results for
coactions of $\hat\sv$:  for example, Equation~\eqref{xr-eq2} 
above can be derived
from Equation~\eqref{xr-eq1} by replacing $V$ by $\hat V$.

Let $V$ be a regular multiplicative unitary, and let $\psi\colon A\to
M(B)$ be a nondegenerate homomorphism which is {\em equivariant\/} for
coactions $\delta_A$ and $\delta_B$ of $\sv$; i.e. such that
$$\delta_B\circ\psi = (\psi\otimes\id)\circ\delta_A.$$
If $(j_B,v_B)$ is the universal covariant pair for $(B,\sv,\delta_B)$
on $B\times\hat\sv$, then $(j_B\circ\psi,v_B)$ is a covariant pair for
$(A,\sv,\delta_A)$, so
we get a nondegenerate homomorphism
$\psi\times\hat\sv = (j_B\circ \psi)\times\mu_B
\colon A\times\hat\sv\to M(B\times\hat\sv)$,
where $\mu_B\colon (\hat\sv)_p\to M(A\times\hat\sv)$ corresponds to
$v_B$ as in \cite[Lemma~2.6]{NgCC1}.
If $\psi$ is equivariant for coactions of $\hat\sv$ on $A$ and $B$, we
likewise get a nondegenerate homomorphism $\psi\times\sv\colon
A\times\sv\to M(B\times\sv)$. 

The analogous result for reduced crossed products, which we will need
in order to define the reduced right-Hilbert bimodule crossed products
in the next section, requires a bit more
work:  

\begin{lem}\label{psi-xr-lem}
Let $V$ be a regular multiplicative unitary on a Hilbert space $H$, and
let $\delta_A$ and $\delta_B$ be coactions of $\sv$ on $C^*$-algebras
$A$ and $B$.  Suppose also that
$\psi\colon A\to M(B)$ is a $\delta_A$ -- $\delta_B$
equivariant nondegenerate homomorphism.  Then there exists a
nondegenerate homomorphism $\psi\times_r\hat\sv\colon
A\times_{\delta_A,r}\hat\sv\to M(B\times_{\delta_B,r}\hat\sv)$ such that
\begin{equation}\label{psi-xr-eq1}
(\psi\times_r\hat\sv)\bigl(\pi_L^A(a)(1\otimes\rho(y))\bigr) =
\pi_L^B(\psi(a))(1\otimes\rho(y))
\end{equation}
for $a\in A$ and $y\in \hat\sv$.  
\end{lem}

\begin{proof}
As in the proof of \cite[Th\'eor\`eme~7.5]{BS-UM}, $A\times_r\hat\sv$
acts nondegenerately on $A\otimes H$, and therefore on $(A\otimes
H)\otimes_AB$, where
$B$ is the standard
right-Hilbert $A$ -- $B$ bimodule arising from $\psi$.
It is straightforward to check that the map $\Phi\colon (A\otimes
H)\otimes_A B \to B\otimes H$ determined by 
$$\Phi((a\otimes\xi)\otimes_A b) = \psi(a)b\otimes\xi$$
is a Hilbert $B$-module isomorphism; thus 
$\L_B((A\otimes H)\otimes_AB)\cong \L_B(B\otimes H)$, so
we obtain
a nondegenerate homomorphism
$\psi\times_r\hat\sv\colon A\times_r\hat\sv\to \L_B(B\otimes H)$
characterized by
$$(\psi\times_r\hat\sv)\bigl(\pi_L^A(a)(1\otimes\rho(y))\bigr)
\bigl(\Phi((c\otimes\xi)\otimes_Ab)\bigr) =
\Phi\bigl(\pi_L^A(a)(1\otimes\rho(y))((c\otimes\xi)\otimes_Ab)\bigr)
$$
for $a,c\in A$, $y\in\hat\sv$, $\xi\in H$, and $b\in B$.  

Now for $a,c,\xi$, and $b$ as above, factor $\xi =
L(x)\eta$ for some $x\in\sv$ and $\eta\in H$, 
and choose $a_i\in A$ and $x_i\in
\sv$ such that $\delta_A(a)(1\otimes x) \approx \sum_i^n a_i\otimes
x_i$.  
Then we have
\begin{eqnarray*}
\lefteqn{
(\psi\times_r\hat\sv)(\pi_L^A(a))\bigl(\Phi((c\otimes\xi)\otimes_Ab)\bigr)}\\
& = & \Phi\bigl(\pi_L^A(a)(c\otimes\xi)\otimes_Ab\bigr)\\
& = & \Phi\bigl((\id\otimes L)\circ\delta_A(a)(c\otimes L(x)\eta)
\otimes_Ab\bigr)\\
& = & \Phi\bigl((\id\otimes L)(\delta_A(a)(1\otimes x))(c\otimes\eta)
\otimes_A b\bigr)\\
& \approx & \sum_i^n \Phi\bigl((\id\otimes L)(a_i\otimes x_i)(c\otimes\eta)
\otimes_A b\bigr)\\
& = & \sum_i^n \Phi\bigl((a_i\otimes L(x_i))(c\otimes\eta)
\otimes_A b\bigr)\\
& = & \sum_i^n \Phi\bigl((a_ic\otimes L(x_i)\eta)\otimes_A b\bigr)\\
& = & \sum_i^n \psi(a_ic)b\otimes L(x_i)\eta\\
& = & \sum_i^n (\psi(a_i)\otimes L(x_i))(\psi(c)b\otimes\eta)\\
& = & (\psi\otimes L)\Bigl(\sum_i^n a_i\otimes x_i\Bigr)
(\psi(c)b\otimes\eta)\\
& \approx & (\psi\otimes L)\bigl(\delta_A(a)(1\otimes x)\bigr)
(\psi(c)b\otimes\eta)\\
& = & (\psi\otimes L)(\delta_A(a))(1\otimes
L(x))(\psi(c)b\otimes\eta)\\
& = & (\psi\otimes L)\circ\delta_A(a)(\psi(c)b\otimes \xi)\\
& = & (\id\otimes L)\circ\delta_B(\psi(a))
\Phi((c\otimes\xi)\otimes_Bb)\\
& = & \pi_L^B(\psi(a))\bigl(\Phi((c\otimes\xi)\otimes_Bb)\bigr),
\end{eqnarray*}
so that $(\psi\times_r\hat\sv)(\pi_L^A(a))
= \pi_L^B(\psi(a))$.  
Since it is straightforward to check that
$(\psi\times_r\hat\sv)(1\otimes\rho(y))
= 1\otimes\rho(y)$ for $y\in
\hat\sv$, this shows that $\psi\times_r\hat\sv$ maps 
$A\times_r\hat\sv$ into
$M(B\times_r\hat\sv)\subseteq \L_B(B\otimes H)$  and also 
establishes Equation~\eqref{psi-xr-eq1}, which in turn makes it evident
that $\psi\times_r\hat\sv$ is nondegenerate.
\end{proof}


\section{Coactions on right-Hilbert bimodules}
\label{rh-bimod-sec}

Let $V$ be a regular multiplicative unitary.  For simplicity, we'll
just write $S$ for $S_V$ and $\hat S$ for $\hat S_V$.
We define
a {\em coaction\/} of the Hopf $C^*$-algebra $S$
on a right-Hilbert $A$ -- $B$ bimodule $X$ to be an imprimitivity
bimodule coaction $(\delta_K,\delta_X,\delta_B)$ of $S$ on $_KX_B$
(\cite[Definition~3.3(a)]{NgCC2}) together with a $C^*$-coaction
$\delta_A$ of $S$ on $A$ such that the associated homomorphism
$\psi\colon A\to M(K)$ is $\delta_A$ -- $\delta_K$ equivariant.  We say
that a right-Hilbert bimodule coaction $(\delta_A,\delta_X,\delta_B)$
is {\em injective\/} if $\delta_A$ and $\delta_B$ are (in which case
$\delta_X$ will be also),
and 
we say it is {\em nondegenerate\/} if $\delta_A$
and $\delta_B$ are nondegenerate
$C^*$-coactions.  

Given an imprimitivity bimodule coaction $(\delta_K,\delta_X,\delta_B)$
of $S$ on $_KX_B$, the rule 
$$\delta_L\twoby{k}{x}{\tilde y}{b} =
\twoby{\delta_K(k)}{\delta_X(x)}{\delta_X(y)\tilde{}}{\delta_B(b)}$$
defines a coaction $\delta_L$ of $S$ on the linking algebra
$L(X)=\ttwoby{K}{X}{\tilde X}{B}$ (\cite[Lemma~3.7]{NgCC2}).  Departing
slightly from Ng, we will define the {\em imprimitivity bimodule
crossed product\/} $_KX_B\times_{\delta_X}\hat S$ to be the corner $j_L(p)
(L(X)\times_{\delta_L}\hat S) j_L(q)$, where $p=\ttwoby{1}{0}{0}{0}$ and
$q=\ttwoby{0}{0}{0}{1}$ are the canonical projections in $M(L(X))$.  
By \cite[Theorem~3.11]{NgCC2}, $_KX_B\times_{\delta_X}\hat S$ is then 
a $K\times_{\delta_K}\hat S$ -- $B\times_{\delta_B}\hat S$ imprimitivity
bimodule which is an imprimitivity bimodule crossed product in Ng's
sense (\cite[Definition~3.5(b)]{NgCC2}), and we have
$$L(X)\times_{\delta_L}\hat S 
\cong \twoby{K\times_{\delta_K}\hat S}{X\times_{\delta_X}\hat
S}{(X\times_{\delta_X}\hat S)\tilde{}}{B\times_{\delta_B}\hat S}
=L(X\times_{\delta_X}\hat S).
$$
Similarly (this time in keeping with Ng), we define the {\em reduced
crossed product\/} $_KX_B\times_{\delta_X,r}\hat S$ to be $(\id\otimes
L)\circ\delta_L(p)(L(X)\times_{\delta_L,r}\hat S)(\id\otimes
L)\circ\delta_L(q)$ (\cite[Remark~3.20(a)]{NgCC2}).  By the proof of
\cite[Proposition~3.19]{NgCC2}, it is a $K\times_{\delta_K,r}\hat
S$ -- $B\times_{\delta_B,r}\hat S$ imprimitivity bimodule, and 
$$L(X)\times_{\delta_L,r}\hat S 
\cong \twoby{K\times_{\delta_K,r}\hat S}{X\times_{\delta_X,r}\hat
S}{(X\times_{\delta_X,r}\hat S)\tilde{}}{B\times_{\delta_B,r}\hat S}
=L(X\times_{\delta_X,r}\hat S).$$

Now given a right-Hilbert bimodule coaction $(\delta_A,\delta_X,\delta_B)$
of $S$ on $_AX_B$, 
the nondegenerate homomorphism $\psi\times\hat S\colon A\times\hat
S\to M(K\times\hat S)$ makes $_KX_B\times\hat S$ 
into a right-Hilbert $A\times\hat
S$ -- $B\times\hat S$ bimodule, which we denote by $_AX_B\times\hat S$
and call the {\em right-Hilbert bimodule crossed product\/} of $_AX_B$
by $S$.  
Similarly, the nondegenerate
homomorphism $\psi\times_r\hat S\colon A\times_r\hat S\to
M(K\times_r\hat S)$ of \lemref{psi-xr-lem}
makes $_KX_B\times_r\hat S$ into 
a right-Hilbert $A\times_r\hat S$ -- $B\times_r\hat S$ bimodule, which we
denote $_AX_B\times_r\hat S$.  
If $V$ comes from a Kac system, we define
right-Hilbert bimodule coactions of $\hat\sv\cong S_{\hat V}$ 
and the right-Hilbert bimodule crossed products $_AX_B\times\sv$ and
$_AX_B\times_r\sv$ by replacing $V$ with $\hat V$ in the above
definitions.  

If $(\delta_K,\delta_X,\delta_B)$ is an
imprimitivity bimodule coaction
of $ S$ on $_KX_B$, it is straightforward to check that the dual
coaction $\hat\delta_L$
of $\hat S_p$ on $L(X)\times\hat S$
restricts to the dual coactions $\hat\delta_K$ and $\hat\delta_B$ on
the diagonal corners 
$K\times\hat S$ and $B\times\hat S$.  The restriction of
$\hat\delta_L$ to the upper right corner
$_KX_B\times\hat S$ gives a map $\hat\delta_X$ such
that $(\hat\delta_K,\hat\delta_X,\hat\delta_B)$ is an imprimitivity
bimodule coaction of $\hat S_p$ on $_KX_B\times\hat S$ which we
call the {\em dual imprimitivity bimodule coaction\/}.  (One can show
that this definition agrees with that given for Hilbert modules in
\cite[Remark~2.18]{NgCC2}.)  
If $(\delta_A,\delta_X,\delta_B)$ is a right-Hilbert bimodule coaction
of $ S$ on
$_AX_B$, we define the {\em dual coaction\/} of $\hat S_p$
on $_AX_B\times\hat S$ to be the dual imprimitivity bimodule coaction
$(\hat\delta_K,\hat\delta_X,\hat\delta_B)$, 
together with the dual $C^*$-coaction
$\hat\delta_A$.  Since
\begin{eqnarray*}
\lefteqn{((\psi\times\hat S)\otimes\id)\circ\hat\delta_A(
j_A(a)\mu_A(y))}\\
& = & ((\psi\times\hat S)\otimes\id)\bigl((j_A(a)\otimes
1)(\mu_A\otimes\id)
(\hat\delta_{V}(y))\bigr)\\
& = & (j_K(\psi(a))\otimes1)(\mu_K\otimes\id)(\hat\delta_{V}(y))\\
& = & \hat\delta_K(j_K(\psi(a))\mu_K(y))\\
& = & \hat\delta_K\circ(\psi\times\hat S)(j_A(a)\mu_A(y))
\end{eqnarray*}
for all $a\in A$, $y\in \hat S_p$, the nondegenerate homomorphism
$\psi\times\hat S\colon A\times\hat S\to M(K\times\hat S)$ is
$\hat\delta_A$ -- $\hat\delta_K$ equivariant, so this is indeed a
right-Hilbert bimodule coaction.  

Given a right-Hilbert bimodule coaction $(\delta_A,\delta_K,\delta_K)$
of $S$ on a standard bimodule $_AK_K$, we have potentially two
different right-Hilbert $A\times\hat S$ -- $K\times\hat S$ bimodules: the
bimodule crossed product $_AK_K\times\hat S$ and the standard bimodule
formed from the $C^*$-algebra crossed product $K\times\hat S$ and
the nondegenerate homomorphism 
$\psi\times\hat S\colon A\times\hat S\to M(K\times\hat S)$.  
The following lemma shows that these coincide; in
other words, a
crossed product of a standard bimodule is a standard
bimodule.

\begin{lem}\label{std-lem}
Let $V$ be a regular multiplicative unitary, and 
let $(\delta_A,\delta_K,\delta_K)$ be a 
right-Hilbert bimodule coaction of $S=S_V$ on a standard bimodule $_AK_K$.  
Then the right-Hilbert bimodule crossed product
$_AK_K\times_{\delta_K}\hat S$ is isomorphic to the
$C^*$-crossed product
$K\times_{\delta_K}\hat S$ as a
right-Hilbert 
$A\times_{\delta_A}\hat S$ -- $K\times_{\delta_K}\hat S$
bimodule.  An analogous statement also holds for the
reduced crossed products. 
\end{lem}

\begin{proof}
Since $_AK_K\times\hat S$ is by definition
$_KK_K\times\hat S$ with the same homomorphism
$\psi\times\hat S$,  it suffices to show that $_KK_K\times\hat S$ is
isomorphic to the $C^*$-crossed product
$K\times_{\delta_K}\hat S$ as a
$K\times_{\delta_K}\hat S$ -- $K\times_{\delta_K}\hat S$
imprimitivity bimodule.
Let $L=\ttwoby{K}{K}{K}{K}$ be the linking algebra for $_KK_K$, and
let 
$\delta_L=\ttwoby{\delta_K}{\delta_K}{\delta_K}{\delta_K}$ be the
associated coaction.  Then 
by \cite[Theorem~3.11]{NgCC2}, 
\begin{equation}\label{std-eq1}
L\times_{\delta_L}\hat S\cong
\twoby{K\times\hat S}{_KK_K\times\hat S}{_KK_K\times
\hat S}{K\times\hat S},
\end{equation}
and the projection $j\ttwoby{1_K}{0}{0}{0}$ in
$M(L\times_{\delta_L}\hat S)$ corresponds to
$\ttwoby{1_{K\times\hat S}}{0}{0}{0}$ under this isomorphism. 

Let $M_2$ denote the $C^*$-algebra of two-by-two matrices over $\bbC$.
Then the canonical isomorphism $\Phi\colon M_2\otimes K\to L$ determined
by $\Phi\left(\ttwoby{a}{b}{c}{d}\otimes k\right) =
\ttwoby{ak}{bk}{ck}{dk}$ is clearly $\id\otimes\delta_K$ -- $\delta_K$
equivariant; thus
\begin{equation}\label{std-eq2}
L\times_{\delta_L}\hat S\cong (M_2\otimes
K)\times_{\id\otimes\delta_K}\hat S.
\end{equation}
Note that this isomorphism takes $j\ttwoby{1_K}{0}{0}{0}\in
M(L\times_{\delta_L}\hat S)$ to $j\left( \ttwoby{1}{0}{0}{0}\otimes
1_K\right)
\in M((M_2\otimes
K)\times_{\id\otimes\delta_K}\hat S)$.  

We next claim that 
\begin{equation}\label{std-eq3}
(M_2\otimes K)\times_{\id\otimes\delta_K}\hat S \cong M_2\otimes
(K\times_{\delta_K}\hat S).
\end{equation}
For if $\iota$ denotes the trivial coaction of $\bbC$ on $M_2$, then
\cite[Proposition~3.2]{NgCC1}\ implies that 
$$(M_2\otimes K)\times_{\iota\otimes\delta_K}\widehat{\bbC\otimes S}
\cong (M_2\times_\iota\hat\bbC)\otimes(K\times_{\delta_K}\hat S),$$
and the right and left sides of this equation are naturally isomorphic
to the right and left sides, respectively, of Equation~\eqref{std-eq3}.  
Under this isomorphism, the projection
$j(\ttwoby{1}{0}{0}{0}\otimes 1_K)\in M((M_2\otimes
K)\times_{\id\otimes\delta_K}\hat S)$ is carried to
$\ttwoby{1}{0}{0}{0}\otimes 1_{K\times\hat S}\in M(M_2\otimes
(K\times_{\delta_K}\hat S))$.  

Combining Equations~(\ref{std-eq1}), (\ref{std-eq2}), 
and (\ref{std-eq3}), we have
$$\twoby{K\times\hat S}{_KK_K\times\hat S}{_KK_K\times
\hat S}{K\times\hat S}\cong 
M_2\otimes (K\times_{\delta_K}\hat S),$$
and since $j\ttwoby{1_K}{0}{0}{0}$ maps to 
$\ttwoby{1}{0}{0}{0}\otimes 1$, it follows that the 
corners $_KK_K\times_{\delta_K}\hat S$ and $K\times_{\delta_K}\hat S$
are isomorphic as
$K\times_{\delta_K}\hat S$ -- $K\times_{\delta_K}\hat S$ imprimitivity
bimodules. 

For the reduced crossed products, 
it again suffices to show that $_KK_K\times_r\hat S$ is isomorphic to
$K\times_r\hat S$ as a $K\times_r\hat S$ -- $K\times_r\hat S$
imprimitivity bimodule.  
By \cite[Remark~3.20(a)]{NgCC2}\ we have
$$L\times_{\delta_L,r}\hat S\cong 
\twoby{K\times_r\hat S}{_KK_K\times_r\hat S}{_KK_K\times_r
\hat S}{K\times_r\hat S},$$
and it follows from equivariance of $\Phi\colon M_2\otimes K\to L$
that
$$L\times_{\delta_L,r}\hat S\cong (M_2\otimes
K)\times_{\id\otimes\delta_K,r}\hat S.$$
Applying \cite[Proposition~3.3]{NgCC1},
which is the reduced version of
\cite[Proposition~3.2]{NgCC1},
we get
$$(M_2\otimes K)\times_{\iota\otimes\delta_K,r}\widehat{\bbC\otimes S}
\cong (M_2\times_{\iota,r}\hat\bbC)\otimes(K\times_{\delta_K,r}\hat S),$$
and hence 
$$(M_2\otimes K)\times_r\hat S \cong M_2\otimes (K\times_r\hat S).$$
Combining these isomorphisms and matching up the projections as
above, it follows that
$_KK_K\times_{\delta_K,r}\hat S\cong K\times_{\delta_K,r}\hat S$.
\end{proof}

For the proof of our main result (\thmref{main}) we will need to know
that the decomposition ${}_AX_B\cong{}_AK\otimes_KX_B$ is equivariant:

\begin{lem}\label{tensor}
Let $V$ be a regular multiplicative unitary, let 
$(\delta_A,\delta_X,\delta_B)$ be a right-Hilbert bimodule
coaction of $S=S_V$ on $_AX_B$, and let $\delta_K$ be the
associated coaction on $K=\K_B(X)$.  Then there exist right-Hilbert
bimodule isomorphisms
$$_AX_B\times_{\delta_X}\hat S \cong (_AK_K\times_{\delta_K}\hat
S)\otimes_{K\times\hat S}(_KX_B\times_{\delta_X}\hat S)$$
and
$$_AX_B\times_{\delta_X,r}\hat S \cong (_AK_K\times_{\delta_K,r}\hat
S)\otimes_{K\times_r\hat S}(_KX_B\times_{\delta_X,_r}\hat S).$$
\end{lem}

\begin{proof}
By definition, $_AX_B\times_{\delta_X}\hat S$ is the imprimitivity
bimodule $_KX_B\times_{\delta_X}\hat S$ with the nondegenerate
homomorphism 
$\psi\times\hat S\colon A\times\hat S\to M(K\times\hat S)$ 
arising
from $\psi\colon A\to M(K)$.  Since $_AK_K\times\hat S$ is
$_KK_K\times\hat S$ with the same map, for the first isomorphism it
suffices to show that $_KX_B\times_{\delta_X}\hat S \cong
(_KK_K\times_{\delta_K}\hat
S)\otimes_{K\times\hat S}(_KX_B\times_{\delta_X}\hat S)$ as
imprimitivity bimodules.  But by \lemref{std-lem},
$_KK_K\times\hat S\cong K\times\hat S$, so the result follows
from the usual cancelation $C\otimes_CY\cong Y$. 

The assertion about the reduced crossed products follows similarly
from \lemref{std-lem}.   
\end{proof}


\section{The fixed-point theorem}
\label{fixed-point}

Based upon the familiar results for actions of compact groups and 
coactions of discrete groups, one 
would guess that the crossed product by
a coaction of a Hopf $C^*$-algebra of compact type is Morita equivalent to the
fixed-point algebra. In \cite{NgME} Ng proves a version of this fixed-point
theorem, and this is crucial for his imprimitivity theorem, which we study in
the next section. 
(We should point out
that the imprimitivity theorem naturally involves a
coaction of $\hat\su$ for a multiplicative unitary $U$ coming from a
{\em discrete\/} Kac system, but is proved by applying 
the fixed-point theorem 
to the corresponding coaction of $S_{\hat U}$, where $\hat U$ is
compact.)
Here we recall Ng's
fixed-point result and establish some relations to
multipliers and bimodules, in preparation for our work with Ng's imprimitivity
theorem in \secref{imprimitivity}.

Let $V$ be a regular multiplicative unitary of compact type such that
$\sv$ has a faithful Haar state $\varphi$. Again 
we'll just write $S$ for $\sv$ and $\hat S$ for
$\hat \sv$.  Let $\delta$ be a coaction of $S$ on a $C^*$-algebra $A$
which is \emph{effective} in the sense that
\[
\clsp\{\delta(A)(A\otimes 1)\}=A\otimes S.
\]
Ng shows in two steps (\cite[Theorem 2.7 and Proposition 2.9]{NgME}) that
the reduced crossed product $A\times_{\delta,r}\hat S$ is Morita
equivalent to the fixed-point algebra $A^\delta$. Since it will
simplify our computations with the imprimitivity bimodule, we will
combine Ng's two  steps into one. 

Ng's strategy is to use a nonunital version of Watatani's $C^*$-basic
construction (\cite{WatIC}). The map $E=E_A=(\id\otimes\varphi)\circ\delta$
is a conditional expectation of $A$ onto $A^\delta$, and so $A$ becomes
a full pre-Hilbert $A^\delta$-module under right multiplication and the
pre-inner product
\[
\rinner{a,b}{A^\delta}=E(a^*b).
\]
The Hausdorff completion of the pre-Hilbert $A^\delta$-module $A$ is a
full Hilbert $A^\delta$-module, denoted $\FF=\FF(A)$. Let $\eta=\eta_A$
be the canonical map of $A$ into $\FF$, and define 
$e_A\in\L_{A^\delta}(\FF)$
and $\lambda=\lambda_A\colon A\to\L_{A^\delta}(\FF)$ by
\[
e_A\eta(a)=\eta(E(a))\midtext{and}\lambda(a)\eta(b)=\eta(ab).
\]
Then the \emph{$C^*$-basic construction} is defined to be the closed span
in $\L_{A^\delta}(\FF)$ of $\lambda(A)e_A\lambda(A)$, and is denoted
$C^*\inner{A,e_A}$. Since
\[
e_A\lambda(a)e_A=\lambda(E(a))e_A\midtext{and}E(a^*)=E(a)^*,
\]
$C^*\inner{A,e_A}$ is a $C^*$-algebra; in fact, a routine
computation shows $C^*\inner{A,e_A}$ coincides with the imprimitivity
algebra $\K_{A^\delta}(\FF)$.  Moreover, a short computation shows that
the left inner product is given on the generators by
\[
\linner{\eta(a),\eta(b)}{C^*\inner{A,e_A}}=\lambda(a)e_A\lambda(b^*).
\]
Therefore, the Hausdorff completion $\FF$ of the
\imp{\spn\{\lambda(A)e_A\lambda(A)\}}{A^\delta} pre-imprimitivity
bimodule $A$ is a \imp{C^*\inner{A,e_A}}{A^\delta} imprimitivity
bimodule.

Ng's first step is to temporarily assume the coaction $\delta$ is
injective. Then the conditional expectation $E$ is faithful, and Ng
proves (\cite[Theorem 2.7]{NgME}) that in this case the map
\[
\lambda(a)e_A\lambda(b)\mapsto\delta(a)(1\otimes\rho(\varphi))\delta(b)
\]
extends to an isomorphism of the $C^*$-basic construction
$C^*\inner{A,e_A}$ onto the reduced crossed product
$A\times_{\delta,r}\hat S$, where
\[
\rho(\varphi)=(\id\otimes\varphi)(V),
\]
which, as Ng observes in \cite[proof of Lemma 2.5]{NgME}, is a member
of $\hat S$.

Ng's second step is to remove the injectivity condition on $\delta$ and
note that, if we put $I=\ker\delta$, there is an injective coaction
$\delta'$ on $A/I$ given by $\delta'(q(a))=(q\otimes\id)\circ\delta(a)$,
where $q\colon A\to A/I$ is the quotient map. Then $\delta'$ is also
effective, $q$ maps $A^\delta$ isomorphically onto $(A/I)^{\delta'}$, and
the reduced crossed products $A\times_{\delta,r}\hat S$ 
and $(A/I)\times_{\delta',r}\hat S$
coincide. Ng deduces as a corollary (\cite[Proposition 2.9]{NgME})
that $A\times_{\delta,r}\hat S$ is still Morita equivalent to~$A^\delta$.

To combine Ng's two steps, note that in the second step the imprimitivity
bimodule $\FF(A/I)$ is the completion of $A/I$ with inner product
\[
\rinner{q(a),q(b)}{(A/I)^{\delta'}}=E_{A/I}(q(a)^*q(b))
=q\circ E_A(a^*b)
=q\bigl(\rinner{a,b}{A^\delta}\bigr).
\]
Since $q$ is faithful on the image of $\rinner{\cdot,\cdot}{A^\delta}$,
$\FF(A/I)$ can be identified with $\FF(A)$. More precisely, the map
$\eta_A(a)\mapsto\eta_{A/I}(q(a))$ is well-defined and
extends to an isomorphism $\Phi$ of the Hilbert
$A^\delta$-module $\FF(A)$ onto the Hilbert $(A/I)^{\delta'}$-module
$\FF(A/I)$, with right coefficient map $q|_{A^\delta}$. Moreover, a
short computation shows
\[
\Phi\bigl(\lambda_A(a)e_A\lambda_A(b)\eta_A(c)\bigr)
=\lambda_{A/I}(q(a))e_{A/I}\lambda_{A/I}(q(b))\Phi(\eta_A(c)),
\]
so $\Phi$ is in fact an isomorphism of the
\imp{C^*\inner{A,e_A}}{A^\delta} imprimitivity bimodule $\FF(A)$ onto
the \imp{C^*\inner{A/I,e_{A/I}}}{(A/I)^{\delta'}} imprimitivity bimodule
$\FF(A/I)$, with left coefficient map determined by
\[
\lambda_A(a)e_A\lambda_A(b)\mapsto
\lambda_{A/I}(q(a))e_{A/I}\lambda_{A/I}(q(b)).
\]
Combining with the isomorphism
\[
\lambda_{A/I}(q(a))e_{A/I}\lambda_{A/I}(q(b))\mapsto
\delta'(q(a))(1\otimes\rho(\varphi))\delta'(q(b))
\]
of $C^*\inner{A/I,e_{A/I}}$ onto $(A/I)\times_{\delta',r}\hat S$, and with the
identification of $A\times_{\delta,r}\hat S$ and 
$(A/I)\times_{\delta',r}\hat S$, we get
an isomorphism
\[
\lambda_A(a)e_A\lambda(b)\mapsto\delta(a)(1\otimes\rho(\varphi))\delta(b)
\]
of $C^*\inner{A,e_A}$ onto $A\times_{\delta,r}\hat S$. Putting all this
together, we have a one-step version of Ng's fixed-point theorem
(\cite[Proposition 2.11]{NgME})
--- although we haven't addressed the case of coactions by $S_p$:

\begin{prop}[\cite{NgME}]\label{Ng-prop}
If $V$ is a regular multiplicative unitary of compact type such that
$S=\sv$ has a faithful Haar state $\varphi$, and if $\delta$ is an effective
coaction of $ S$ on $A$, then $A$ is a pre-imprimitivity bimodule
between the pre-$C^*$-algebra
$B=\spn\{\delta(A)(1\otimes\rho(\varphi))\delta(A)\}$
and the fixed-point
algebra $A^\delta$, with operations given for $a,b,c\in A$ and $d\in
A^\delta$ by 
\begin{align*}
\bigl(\delta(a)(1\otimes\rho(\varphi))\delta(b)\bigr)\cdot c&=aE(bc)\\
a\cdot d&=ad\\
\linner{a,b}{B}&=\delta(a)(1\otimes\rho(\varphi))\delta(b^*)\\
\rinner{a,b}{A^\delta}&=E(a^*b).
\end{align*}
Consequently, the Hausdorff completion $\FF(A)$ of $A$ is an
\imp{A\times_{\delta,r}\hat S}{A^\delta} imprimitivity bimodule.
\end{prop}

Now let $(\delta_A,\delta_X,\delta_B)$ be a coaction of $S$ on an
\imp{A}{B} imprimitivity bimodule $X$, let $L=L(X)$ be the linking
algebra, and let $\delta_L$ be the associated coaction of $S$ on $L$.  
Then we have 
\begin{align*}
\delta_L(L)(L\otimes 1)
&\supseteq\begin{pmatrix}\delta_A(A)&\delta_X(X)\\
{\delta_X(X)\tilde{\ }}&\delta_B(B)\end{pmatrix}
\begin{pmatrix}A\otimes 1&0\\0&B\otimes 1\end{pmatrix}\\
&=\begin{pmatrix}\delta_A(A)(A\otimes 1)&\delta_X(X)\cdot(B\otimes 1)\\
{\delta_X(X)\tilde{\ }}\cdot(A\otimes 1)&
\delta_B(B)(B\otimes 1)\end{pmatrix}
\end{align*}
and
\begin{align*}
\delta_X(X)\cdot(B\otimes 1)
&=\delta_X(X\cdot B)\cdot(B\otimes 1)
\\&=\delta_X(X)\cdot\delta_B(B)(B\otimes 1);
\end{align*}
it follows from this (and by symmetry) that $\delta_L$ is effective
whenever $\delta_A$ and $\delta_B$ are.

Let $p=\intextmatrix{1&0\\0&0}\in M(L)$.
We will need the following result in the next section.

\begin{lem}
\label{basic corner}
Let $V$ be a regular multiplicative unitary of compact type such that
$S=\sv$ has a faithful Haar state $\varphi$, and with notation as above,
suppose that $\delta_A$ and $\delta_B$ are effective.  
Then the inclusion $A\hookrightarrow L$ extends to an isomorphism $\Phi$ of
the \imp{A\times_{\delta_A,r}\hat S}{A^{\delta_A}} 
imprimitivity bimodule $\FF(A)$
onto the 
\imp{\delta_L(p)(L\times_{\delta_L,r}\hat S)\delta_L(p)}{pL^{\delta_L}p}
imprimitivity bimodule $\delta_L(p)\cdot\FF(L)\cdot p$.
\end{lem}

\begin{proof}
Let's first make sure we understand all the components of the statement
of the lemma. On the right side of $\FF(L)$ we regard $p$ as an element
of $M(L^{\delta_L})$, which naturally embeds in $M(L)$. Since the projections $p$
in $M(L^{\delta_L})$ and $\delta_L(p)=p\otimes1$ in $M(L\times_r\hat S)$ are full,
$\delta_L(p)\cdot\FF(L)\cdot p$ is indeed a
\imp{\delta_L(p)(L\times_r\hat S)\delta_L(p)}{pL^{\delta_L}p}
imprimitivity
bimodule. We have
\begin{align*}
\delta_L(p)(L\times_r\hat S)\delta_L(p)
&=(p\otimes 1)
\clsp\{\delta_L(L)(1_{M(L)}\otimes\hat S)\}(p\otimes 1)
\\&=\clsp\{(p\otimes 1)\delta_L(L)(p\otimes 1)(1_{M(A)}\otimes\hat S)\}
\\&=\clsp\{\delta_L(p)\delta_L(L)\delta_L(p)(1_{M(A)}\otimes\hat S)\}
\\&=\clsp\{\delta_L(pLp)(1_{M(A)}\otimes\hat S)\}
\\&=\clsp\{\delta_A(A)(1_{M(A)}\otimes\hat S)\}
\\&=A\times_r\hat S,
\end{align*}
where we have used $p=1_{M(A)}$.
On the other side, since $E_L|_A=E_A$ and the natural extension of $E_L$
to $M(L)$ is a conditional expectation onto $M(L^{\delta_L})$, we have
\begin{align*}
pL^{\delta_L}p
=pE_L(L)p
=E_L(pLp)
=E_A(A)
=A^{\delta_A}.
\end{align*}

Thus, it suffices to show the inclusion $A\hookrightarrow L$ respects
the right inner products and the left module multiplications. For the
inner products, if $a,b\in A$ then
\begin{align*}
\rinner{a,b}{L^{\delta_L}}
&=E_L(a^*b)
=E_A(a^*b)
=\rinner{a,b}{A^{\delta_A}}.
\end{align*}
Turning to the left module multiplications, first note that
\[
\delta_L(p)(1_{M(L)}\otimes\rho(\varphi))\delta_L(p)
=p\otimes\rho(\varphi)
=1_{M(A)}\otimes\rho(\varphi),
\]
so for $a,b\in A$ we have
\[
\delta_L(a)(1_{M(L)}\otimes\rho(\varphi))\delta_L(b)
=\delta_A(a)(1_{M(A)}\otimes\rho(\varphi))\delta_A(b).
\]
Hence, for $a,b,c\in A$ we have
\begin{align*}
\bigl(\delta_L(a)(1_{M(L)}\otimes\rho(\varphi))\delta_L(b)\bigr)\cdot c
&=aE_L(bc)
=aE_A(bc)
\\&=\bigl(\delta_A(a)(1_{M(A)}\otimes\rho(\varphi))\delta_A(b)\bigr)\cdot c,
\end{align*}
and we're done.
\end{proof}

We'll also need the following lemma concerning standard bimodules.

\begin{lem}
\label{hom}
Let $V$ be a regular multiplicative unitary of compact type such that
$S=\sv$ has a faithful Haar state $\varphi$.  
If $\psi\colon A\to M(B)$ 
is a nondegenerate homomorphism which is
equivariant for effective
coactions $\delta_A$ and $\delta_B$ of $S$, then $\psi$
extends to a nondegenerate imprimitivity bimodule homomorphism 
$\Psi\colon \FF(A)\to M(\FF(B))$
with coefficient maps $\psi\times_r\hat S$
and $\psi|_{A^{\delta_A}}$.
\end{lem}

\begin{proof}
By \cite[Lemma 5.1]{KQR-DR}, it's enough to show $\psi$
preserves both module multiplications and inner products. For $a,b,c\in
A$ and $d\in A^{\delta_A}$ we have
\begin{align*}
\psi\Bigl(\bigl(\delta_A(a)(1\otimes\rho(\varphi))\delta_A(b)\bigr)
\cdot c\Bigr)
&=\psi(aE_A(bc))
=\psi(a)\psi\circ E_A(bc)
\\&=\psi(a)E_B\circ\psi(bc)
=\psi(a)E_B\bigl(\psi(b)\psi(c)\bigr)
\\&=\bigl(\delta_B(\psi(a))(1\otimes\rho(\varphi))\delta_B(\psi(b))\bigr)
\cdot\psi(c)
\\&=(\psi\times_r\hat S)
\bigl(\delta_A(a)(1\otimes\rho(\varphi))\delta_A(b)\bigr)\cdot\psi(c),
\end{align*}
\[
\psi(a\cdot d)=\psi(ad)=\psi(a)\psi(d)=\psi(a)\cdot\psi(d),
\]
\begin{align*}
\linner{\psi(a),\psi(b)}{M(B\times_r\hat S)}
&=\delta_B(\psi(a))(1\otimes\rho(\varphi))\delta_B(\psi(b))
\\&=(\psi\times_r\hat S)
\bigl(\delta_A(a)(1\otimes\rho(\varphi))\delta_A(b)\bigr)
\\&=(\psi\times_r\hat S)\bigl(\linner{a,b}{A\times_r\hat S}\bigr),
\end{align*}
and
\begin{align*}
\rinner{\psi(a),\psi(b)}{M(B^{\delta_B})}
&=E_B\bigl(\psi(a)^*\psi(b)\bigr)
=E_B\circ\psi(a^*b)
\\&=\psi\circ E_A(a^*b)
=\psi|_{A^{\delta_A}}\bigl(\rinner{a,b}{A^{\delta_A}}\bigr).
\end{align*}
\end{proof}


\section{Equivariance and imprimitivity}
\label{imprimitivity}

We now turn to the result we call Ng's imprimitivity theorem
(\cite[Theorem~3.4]{NgME}).  This is an analogue, for multiplicative
unitaries of \emph{discrete} type (in fact, coming from discrete Kac
systems), of Green's and Mansfield's imprimitivity theorems for group
actions and coactions, respectively.  
Our main theorem (\thmref{main}) 
says that Ng's theorem is compatible with equivariant
right-Hilbert bimodules; 
we begin by
introducing the notation and construction of Ng's imprimitivity
bimodule. 

Let $U$, $V$, and $W$ be
multiplicative unitaries coming from discrete Kac systems, and
assume $W$ is a normal submultiplicative unitary of $V$ and $U$ is the
corresponding quotient (see \cite[Definition~3.2]{NgME});
this implies that there exist surjective Hopf $*$-homomorphisms
$L_{V,W}\colon S_V\to S_W$ and $\rho_{V,U}\colon (\hat\sv)_p\to (\hat
S_U)_p$.  
Thus any coaction $\delta$ of $\sv$ on $A$ can be restricted to a
coaction $\delta| = (\id\otimes L_{V,W})\circ\delta$ of $S_W$ on $A$,
and any dual coaction 
$\hat\delta$ of $(\hat\sv)_p$ on $A\times_\delta \hat\sv$
can be restricted to a coaction $\hat\delta| =
(\id\otimes\rho_{V,U})\circ\hat\delta$ of 
$(\hat S_U)_p$ on $A\times_\delta \hat\sv$.  We can pass to the
corresponding coaction of the reduced $C^*$-algebra $\hat\su$ without
changing either the crossed product or the fixed-point algebra,
and we continue to denote this coaction by $\hat\delta|$.

Now assume further that $W$ is amenable, and 
let $\rho_{W,V}\colon (\hat\sw)_p\to (\hat\sv)_p$ be the Hopf
$*$-homomorphism vouchsafed
by the normality of $W$ in $V$.  
Ng shows that 
if $\delta$ is nondegenerate,
the nondegenerate homomorphism 
$$\phi_A= j_A^V\times(\mu_A^V\circ\rho_{W,V})\colon 
A\times_{\delta|}\hat\sw\to M(A\times_\delta\hat\sv)$$
is actually an isomorphism of $A\times_{\delta|}\hat\sw$ onto the
fixed-point algebra $(A\times_\delta\hat\sv)^{\hat\delta|}$, and that
the restricted dual coaction $\hat\delta|$ of $\hat\su$ is effective.
Viewing this as an effective coaction of $S_{\hat U}$ with $\hat U$
compact, \propref{Ng-prop} provides an
$A\times_\delta\hat\sv\times_{\hat\delta|,r}\su$ --
$(A\times_\delta\hat\sv)^{\hat\delta|}$ imprimitivity bimodule
$\FF(A\times_\delta\hat\sv)$;
using the isomorphism $\phi_A$, this becomes an
$A\times_\delta\hat\sv\times_{\hat\delta|,r}\su$ -- 
$A\times_{\delta|}\hat\sw$
imprimitivity bimodule which we denote by $N(A)$.  

With notation as below, we define the map
$\delta_X|\colon X\to M(X\otimes\sw)$ to
be $(\id\otimes L_{V,W})\circ\delta_X$; it is straightforward to check
that then $(\delta_K|,\delta_X|,\delta_B|)$ is an imprimitivity
bimodule coaction of $\sw$ on $_KX_B$ and that $\psi\colon A\to M(K)$
is $\delta_A|$ -- $\delta_K|$ equivariant.  We call the resulting
right-Hilbert bimodule coaction $(\delta_A|,\delta_X|,\delta_B|)$ of
$\sw$ the {\em restricted coaction\/} from $\sv$.  The restricted dual
right-Hilbert bimodule coaction
$(\hat\delta_A|,\hat\delta_X|,\hat\delta_B|)$ of $(\hat\su)_p$ is defined
similarly; its reduction to $\hat\su$ is also denoted
$(\hat\delta_A|,\hat\delta_X|,\hat\delta_B|)$.  

\begin{thm}
\label{main}
Let $U$, $V$, and $W$ be multiplicative unitaries coming from discrete
Kac systems, with $W$ an amenable
 normal submultiplicative unitary of $V$ and
$U$ the corresponding quotient. Let  
$(\delta_A,\delta_X,\delta_B)$ be an injective, nondegenerate coaction
of $\sv$ on a right-Hilbert \imp{A}{B} bimodule $X$, and suppose that
the associated coaction $\delta_K$ on the imprimitivity algebra of $X$
is also nondegenerate.  Then the diagram
\begin{equation}
\label{diagram}
\begin{diagram}
\node{A\times_{\delta_A}\hat \sv\times_{\hat\delta_A|,r}\su}
\arrow{e,t}{N(A)}
\arrow{s,l}{_AX_B\times_{\delta_X}\hat \sv\times_{\hat\delta_X|,r}\su}
\node{A\times_{\delta_A|}\hat \sw}
\arrow{s,r}{_AX_B\times_{\delta_X|}\hat \sw}\\
\node{B\times_{\delta_B}\hat \sv\times_{\hat\delta_B|,r}\su}
\arrow{e,b}{N(B)}
\node{B\times_{\delta_B|}\hat \sw}
\end{diagram}
\end{equation}
commutes in the sense that
\[
N(A)\otimes_{A\times\hat \sw}(_AX_B\times\hat \sw)\cong
(_AX_B\times\hat \sv\times_r \su)\otimes_{B\times\hat \sv\times_r \su}N(B)
\]
as right-Hilbert \imp{A\times\hat \sv\times_r \su}{B\times\hat \sw}
bimodules.
\end{thm}

\begin{proof}
By definition we have an imprimitivity bimodule $_KX_B$, a nondegenerate
homomorphism $\psi\colon A\to M(K)$, and a 
coaction $\delta_K$ of $\sv$ on $K=\K_B(X)$ such that
$\psi$ is $\delta_A$ -- $\delta_K$ equivariant and 
$(\delta_K,\delta_X,\delta_B)$ is an imprimitivity bimodule coaction of
$\sv$ on $_KX_B$ which is nondegenerate by assumption.  
Our strategy will be to prove a
version of Diagram~\eqref{diagram}\ for the imprimitivity bimodule $_KX_B$, a
version for the standard bimodule $_AK_K$, and then to combine them
using the decomposition 
\lemref{tensor}.  

First consider the imprimitivity bimodule $_KX_B$:
let $L=L(X)$ be the linking algebra of $X$, let 
$p=\intextmatrix{1&0\\0&0}$ and $q=\intextmatrix{0&0\\0&1}$ be the
canonical projections in $M(L)$, and let $\delta_L$ be the associated
coaction of $\sv$ on $L$. 
Then $\delta_L$ is injective
since $\delta_K$ and $\delta_B$ (hence also $\delta_X$) are. 
Since $\delta_B$ is nondegenerate, we have
\begin{eqnarray*}
\clsp\{\delta_X(X)\cdot(1\otimes S)\}
& = & \clsp\{\delta_X(X\cdot B)\cdot(1\otimes S)\}\\
& = & \clsp\{\delta_X(X)\cdot\delta_B(B)(1\otimes S)\}\\
& = & \clsp\{\delta_X(X)\cdot(B\otimes S)\}\\
& = & X\otimes S;
\end{eqnarray*}
similarly, the nondegeneracy of $\delta_K$ implies that
$\clsp\{\widetilde{\delta_X(X)}\cdot(1\otimes S)\} 
= \widetilde{X\otimes S}$.  It follows
easily that $\delta_L$ is nondegenerate as well. 
Thus, by \cite[Theorem~3.4]{NgME}, we have a \imp{K\times_{\delta_K}\hat
\sv\times_{\hat\delta_K|,r}\su}{K\times_{\delta_K|}\hat \sw}
imprimitivity bimodule $N(K)$, and an \imp{L\times_{\delta_L}\hat
\sv\times_{\hat\delta_L|,r}\su}{L\times_{\delta_L|}\hat \sw}
imprimitivity bimodule $N(L)$. 
We claim that
\begin{equation}\label{main-eq1}
N(K)\otimes_{K\times\hat\sw}(_KX_B\times\hat\sw)\cong
(_KX_B\times\hat\sv\times_r\su)\otimes_{B\times\hat\sv\times_r\su} N(B)
\end{equation}
as $K\times\hat S_V\times_r S_U$ -- $B\times\hat S_W$
imprimitivity bimodules.  

Now $\delta_L|=
\ttwoby{\delta_K|}{\delta_X|}{\delta_{\tilde X}|}{\delta_B|}$, so 
\begin{equation}\label{main-eq3}
L(X)\times_{\delta_L|}\hat\sw\cong L(X\times_{\delta_X|}\hat\sw).
\end{equation}
Also, $(L(X)\times\hat\sv,(\hat\sv)_p,\hat\delta_L)\cong
(L(X\times\hat\sv),(\hat\sv)_p,\epsilon_L)$, where
$\epsilon_L=\ttwoby{\hat\delta_K}{\hat\delta_X}{\hat\delta_{\tilde
X}}{\hat\delta_B}$.  It follows that the coactions $\hat\delta_L|$ and
$\epsilon_L| = \ttwoby{\hat\delta_K|}{\hat\delta_X|}{\hat\delta_{\tilde
X|}}{\hat\delta_B|}$ of $(\hat\su)_p$ are isomorphic, and therefore that
their reductions are, so that 
\begin{equation}\label{main-eq4}
L(X)\times_{\delta_L}\hat\sv\times_{\hat\delta_L|,r}\su
\cong L(X\times_{\delta_X}\hat\sv)\times_{\epsilon_L|,r}\su
\cong L(X\times_{\delta_X}\hat\sv\times_{\hat\delta_X|,r}\su).
\end{equation}
Let $p_W,q_W\in M(L(X\times\hat\sw))$ and $p_U,q_U\in
M(L(X\times\hat\sv\times_r\su))$ be the canonical projections.  Then
using Equations~(\ref{main-eq3}) and~(\ref{main-eq4}) to view $N(L)$ as an 
$L(X\times\hat\sv\times_r\su)$ -- $L(X\times\hat\sw)$ imprimitivity
bimodule, \cite[Lemma~4.6]{ER-ST} gives us a
$K\times\hat\sv\times_r\su$ -- $B\times\hat\sw$ imprimitivity bimodule
isomorphism 
\begin{align*}
(p_U\cdot N(L)\cdot p_W)&\otimes_{K\times\hat\sw}(_KX_B\times\hat\sw)\\
&\cong (_KX_B\times\hat\sv\times_r\su)\otimes_{B\times\hat\sv\times_r\su} 
(q_U\cdot N(L)\cdot q_W).
\end{align*}

Thus, in order to establish Equation~\eqref{main-eq1} we only need 
imprimitivity bimodule isomorphisms
$ p_U\cdot N(L)\cdot p_W\cong N(K)$ and 
$q_U\cdot N(L)\cdot q_W\cong N(B)$,
and by symmetry it suffices to prove the first.
Now the isomorphism $L(X\times\hat\sv\times_r\su)\cong
L(X\times\hat\sv)\times_r\su$ takes $p_U$ to $\epsilon_L(p_V)$, and the
isomorphisms $L(X\times\hat\sw)\cong L(X)\times\hat\sw\cong
(L(X)\times\hat\sv)^{\hat\delta_L|}\cong L(X\times\hat\sv)^{\epsilon_L|}$
carry $p_W$ to $p_V$.  
Therefore, \lemref{basic corner} (applied to the coaction of $S_{\hat U}$
on $L\times_r\hat\sv$ equivalent to $\hat\delta_L|$)
tells us that
\[
p_U\cdot N(L)\cdot p_W
\cong\epsilon_L(p_V)\cdot\FF(L\times_r\hat \sv)\cdot p_V
\cong \FF(K\times_r\hat \sv)
\cong N(K),
\]
which gives Equation~\eqref{main-eq1}. 

Next we consider the standard bimodule $_AK_K$ with the right-Hilbert
bimodule coaction $(\delta_A,\delta_K,\delta_K)$.  
We claim that 
\begin{equation}\label{main-eq2}
N(A)\otimes_{A\times\hat \sw}(_AK_K\times\hat \sw)\cong
(_AK_K\times\hat \sv\times_r \su)\otimes_{K\times\hat \sv\times_r \su}N(K)
\end{equation}
as right-Hilbert \imp{A\times\hat \sv\times_r \su}{K\times\hat \sw}
bimodules;
by \cite[Lemma 5.3]{KQR-DR} and \lemref{std-lem}, 
it's enough to show that there is a
nondegenerate imprimitivity bimodule homomorphism $\Psi$ from $N(A)$
to $M(N(K))$ with coefficient maps $\psi\times\hat \sv\times_r\su$
and $\psi\times\hat \sw$. 
Applying \lemref{hom} to the nondegenerate
homomorphism $\psi\times\hat\sv\colon
A\times\hat\sv\to M(K\times\hat\sv)$, which is equivariant for the
coactions (of $S_{\hat U}$ equivalent to) 
$\hat\delta_A|$ and $\hat\delta_K|$ of $\hat\su$, we obtain a
nondegenerate imprimitivity bimodule homomorphism 
$\Psi\colon \FF(A\times\hat\sv)\to M(\FF(K\times\hat\sv))$ with
coefficient maps $\psi\times\hat\sv\times_r\su$ and
$(\psi\times\hat\sv)|_{(A\times\hat\sv)^{\hat\delta_A|}}$.
Now by definition, $\psi\times\hat\sv =
(j_K^V\circ\psi)\times\mu_K^V$, and Ng's isomorphism $\phi_A\colon
A\times\hat\sw\to (A\times\hat\sv)^{\hat\delta_A|}$ is
$j_A^V\times(\mu_A^V\circ\rho_{W,V})$.  Thus,
\begin{eqnarray*}
\phi_K\circ(\psi\times\hat\sw)
& = & (j_K^V\times(\mu_K^V\circ\rho_{W,V}))\circ 
((j_K^W\circ\psi)\times\mu_K^W)\\
& = & (j_K^V\circ\psi)\times(\mu_K^V\circ\rho_{W,V})\\
& = & \bigl((j_K^V\circ\psi)\times\mu_K^V\bigr)\times
\bigl(j_A^V\times(\mu_A^V\circ\rho_{W,V})\bigr)\\
& = & (\psi\times\hat\sv)\circ\phi_A.
\end{eqnarray*}
This shows that the isomorphisms $\phi_A\colon A\times\hat\sw\to
(A\times\hat\sv)^{\hat\delta_A|}$ and $\phi_K\colon K\times\hat\sw\to
(K\times\hat\sv)^{\hat\delta_K|}$ carry
the coefficient map
$(\psi\times\hat\sv)|_{(A\times\hat\sv)^{\hat\delta_A|}}$ to
$\psi\times\hat\sw$; in
other words, viewed as a map of $N(A)$  into $M(N(K))$, $\Psi$ is a
nondegenerate imprimitivity bimodule homomorphism with coefficient maps
$\psi\times\hat\sv\times_r\su$ and $\psi\times\hat\sw$, which establishes
Equation~\eqref{main-eq2}.

We now have a prism
\begin{equation*}
\renewcommand{\dgeverynode}{\scriptstyle}
\renewcommand{\dgeverylabel}{\scriptscriptstyle}
\dgARROWLENGTH=0em
\begin{diagram}
\node{A\times\hat\sv\times_r\su}
	\arrow[3]{e,t}{N(A)}
	\arrow{sse,t}{_AK_K\times\hat\sv\times_r\su}
	\arrow[4]{s,l}{_AX_B\times\hat\sv\times_r\su}
\node[3]{A\times\hat\sw}
	\arrow{sse,t}{_AK_K\times\hat\sw}
	\arrow[2]{s,-}\\
\node{}\\
\node[2]{K\times\hat\sv\times_r\su}
	\arrow[3]{e,t}{N(K)}
	\arrow{ssw,b}{_KX_B\times\hat\sv\times_r\su}
\node[2]{}
	\arrow[2]{s,l}{_AX_B\times\hat\sw}
\node{K\times\hat\sw}
	\arrow{ssw,b}{_KX_B\times\hat\sw}\\
\node{}\\
\node{B\times\hat\sv\times_r\su}
	\arrow[3]{e,t}{N(B)}
\node[3]{B\times\hat\sw}
\end{diagram}
\end{equation*}
in which the front two faces commute by the above arguments, and the
commutativity of the 
back face is the desired result; it only remains to show that the two
side triangles commute.  That is, we need to know that 
\[
_AX_B\times\hat \sv\times_r \su
\cong(_AK_K\times\hat \sv\times_r \su)
\otimes_{K\times\hat \sv\times_r \su}
(_KX_B\times\hat \sv\times_r \su)
\]
and
\[
_AX_B\times\hat \sw
\cong(_AK_K\times\hat \sw)
\otimes_{K\times\hat \sw}
(_KX_B\times\hat \sw)
\]
as right-Hilbert \imp{A\times\hat \sv\times_r \su}
{B\times\hat \sv\times_r \su} and \imp{A\times\hat \sw}
{B\times\hat \sw} bimodules, respectively. 
But this follows from \lemref{tensor}, in the first case 
applied to the coaction $(\hat\delta_A|,\hat\delta_X|,\hat\delta_B|)$ of
$\hat\su$ on $_AX_B\times\hat\sv$ and then using
$K\times\hat\sv\cong {}_AK_K\times\hat\sv$ from \lemref{std-lem}.
\end{proof}

\begin{rmk}
For an imprimitivity bimodule coaction $(\delta_K, \delta_X,
\delta_B)$ of a Hopf $C^*$-algebra, 
it is probably true that $\delta_K$ is nondegenerate
whenever $\delta_B$ is;  this would simplify the hypotheses of
\thmref{main} somewhat.  Unfortunately, we have been unable to find a
proof.  (For group coactions, it is true --- see
\cite[Proposition~2.3]{KQ-IC} --- and the proof is fairly nontrivial.)
It may be possible to finesse the problem, but since our
main point here is to illustrate our approach to
imprimitivity theorems, we have chosen 
not to get mired in nondegeneracy issues.  
\end{rmk}


\ifx\undefined\bysame
\newcommand{\bysame}{\leavevmode\hbox to3em{\hrulefill}\,}
\fi

\end{document}